\documentclass[12pt,a4paper]{article}
\usepackage[T1]{fontenc}
\usepackage[english]{babel}
\usepackage[centertags]{amsmath}
\usepackage{amsfonts}
\usepackage{amssymb}
\usepackage{amsthm}
\usepackage{newlfont}
\usepackage{graphicx}
\usepackage{color}
\usepackage{float}
\usepackage[inline]{enumitem}
\usepackage{url}
\usepackage{mathtools}

\newtheorem{theorem}{Theorem}[section]
\newtheorem{lemma}[theorem]{Lemma}
\newtheorem{corollary}[theorem]{Corollary}
\newtheorem{proposition}[theorem]{Proposition}
\newtheorem{namedtheorem}{Theorem}

\theoremstyle{definition}
\newtheorem{definition}[theorem]{Definition}
\newtheorem{remark}[theorem]{Remark}
\newtheorem{example}[theorem]{Example}
\newtheorem{question}[theorem]{Question}

\numberwithin{equation}{section}

\newcommand{\blank}{{\mspace{1mu}\cdot\mspace{1mu}}}
\newcommand{\ZZ}{\mathbb{Z}}
\newcommand{\NN}{\mathbb{N}}
\newcommand{\KK}{\mathbb{K}}
\newcommand{\FF}{\mathbb{F}}
\newcommand{\LL}{\mathbb{L}}
\newcommand{\QQ}{\mathbb{Q}}
\newcommand{\RR}{\mathbb{R}}

\newcommand{\pp}{\mathfrak{p}}
\newcommand{\mm}{\mathfrak{m}}
\newcommand{\OO}{\mathcal{O}}
\newcommand{\acts}{{\mspace{1mu}\curvearrowright\mspace{1mu}}}

\DeclareMathOperator{\Iso}{\mathrm{Iso}}

\DeclareMathOperator{\supp}{\mathrm{supp}}

\DeclareMathOperator{\GL}{\mathrm{GL}}
\DeclareMathOperator{\SL}{\mathrm{SL}}
\DeclareMathOperator{\SA}{\mathrm{SA}}

\DeclareMathOperator{\M}{\mathrm{M}}
\DeclareMathOperator{\SO}{\mathrm{SO}}
\DeclareMathOperator{\IA}{\mathrm{Iso_{Aff}}}

\title{Measurable equidecomposability and the Banach--Ruziewicz problem over locally compact non-Archimedean valued fields}

\author{Kamil Orzechowski\\
University of Rzesz\'ow\\ 
Rejtana 16c\\
35-959 Rzesz\'ow, Poland\\
\texttt{kamilo@dokt.ur.edu.pl}\\
ORCID: 0000-0003-0786-1463}

\date{November 24, 2026}

\begin{document}

\maketitle

\begin{abstract}
We prove a measurable analogue of the Banach--Tarski paradox over locally compact non-Archimedean valued fields. Let $\KK$ be the completion of a number field at a non-Archimedean place (e.g., $\KK=\QQ_p$) and $n\ge 2$. Then any two bounded Haar measurable subsets of $\KK^n$ with nonempty interior and equal Haar measure can be cut into finitely many Haar measurable pieces and rearranged into each other by elements of a certain countable group $\Gamma$ of affine isometries. We also show that the Haar measure is, up to a constant factor, the unique nonzero finitely additive $\Gamma$-invariant measure on the bounded Haar measurable subsets of $\KK^n$, which solves a $p$-adic version of the Banach--Ruziewicz problem. Both results follow from the fact that the action of $\Gamma$ on $\KK^n$ is expanding, a property closely related to the spectral gap.
For $n=1$ and arbitrary characteristic, we obtain analogous results when the residue field is not prime, using a countable group of isometries that are not necessarily affine.

  \medskip

  {\footnotesize

  \smallskip

  }
\end{abstract}

\section{Introduction}

Given an action of a group $G$ on a set $X$, two subsets $A,B \subseteq X$ are called (set-theoretically) {\em $G$-equidecomposable} if $A$ can be partitioned into finitely many pieces that can be rearranged, using transformations from $G$, into a partition of $B$. The most famous result concerning equidecomposability is the Banach--Tarski paradox \cite[Theorem~24]{BT}, which states that, for any $n\ge 3$, any two bounded subsets $A,B\subseteq \RR^n$ with nonempty interior are $\Iso(\RR^n)$-equidecomposable, where $\Iso(\RR^n)$ is the group of Euclidean isometries of $\RR^n$.

The result of Banach and Tarski inspired several questions as to whether set-theoretic equidecomposability can be strengthened by imposing additional conditions on the pieces used in the decomposition. Dougherty and Foreman proved \cite{DF} that if such $A,B\subseteq \RR^n$ have additionally the Baire property, they are $\Iso(\RR^n)$-equidecomposable using pieces with the Baire property.
Another question concerns {\em measurable equidecomposability}, that is equidecomposability using Lebesgue measurable pieces. If two measurable subsets $A,B\subseteq \RR^n$ are measurably $\Iso(\RR^n)$-equidecomposable, they clearly have the same Lebesgue measure. Grabowski, M\'ath\'e and Pikhurko proved recently \cite[Theorem~1.2]{GrMaPikhurko} that the converse also holds for $n\ge 3$ and $A,B$ being bounded with nonempty interior.
Moreover, they provided a sufficient criterion for measurable equidecomposability for more general measure-preserving group actions. Their criterion involves the notion of an {\em expanding action}. If a measure-preserving group action by homeomorphisms $G\acts (X,\mu)$ is expanding and {\em paradoxical} (i.e., any two $A,B\subseteq X$ with nonempty interior and compact closure are $G$-equidecomposable), then any two such $A,B$ that are additionally measurable with $\mu(A)=\mu(B)$ are measurably $G$-equidecomposable.

Recently, we studied an analogue of the Banach--Tarski paradox in the new context of normed spaces over non-Archimedean valued fields. In particular, we proved that if $\KK$ is a locally compact non-Archimedean, nontrivially valued field (e.g., $\KK$ is the field $\QQ_p$ of $p$-adic numbers for any prime number $p$) and $n\ge 2$, then any two bounded subsets of $\KK^n$ (equipped with the maximum norm) with nonempty interior are equidecomposable with respect to a certain group of affine isometries of $\KK^n$ \cite[Theorem~5.7]{KO-RACSAM}.

All of this has inspired us to investigate whether this set-theoretical equidecomposability can be upgraded to measurable equidecomposability in the normed spaces $\KK^n$ ($n \ge 2$), equipped with the Haar measure (an analogue---in fact, a generalization---of the Lebesgue measure on $\RR^n$).
In this paper we answer this question affirmatively in the case when $\KK$ has characteristic $0$ (and hence is the completion $K_{\pp}$ of a number field $K$ with respect to a non-Archimedean place $\pp$).

\begin{namedtheorem}\label{thm:A}
Let $\pp$ be a non-Archimedean place of a number field $K$ and $n\ge 2$. Let $\Gamma:=\SL(n,\OO_K)\ltimes K^{n}$, where $\OO_K$ is the ring of integers of $K$. Two Haar measurable subsets $A,B\subseteq K_{\pp}^n$ with compact closure and nonempty interior are Haar-measurably $\Gamma$-equidecomposable if and only if $A$ and $B$ are of equal Haar measure.
\end{namedtheorem}

Corollary~\ref{cor:final}~(i) contains the precise statement for a more general family of Haar measurable sets (not necessarily with nonempty interior, see Example~\ref{ex:empty_interior}), and its version with the Baire property requirement. We emphasize that $\Gamma$ is a countable group of affine isometries of $(K_{\pp}^n,\lVert\blank\rVert_{\infty})$. Theorem~\ref{thm:A} follows from our main technical result---Theorem~\ref{th:main_expanding}---saying that the action of $\Gamma$ on $K_{\pp}^n$, equipped with the Haar measure, is expanding in the sense of Definition~\ref{def:expanding_action}.
Note that if one weakens the notion of measurable equidecomposability by allowing countable partitions, an analogous result is already known; in fact, translations of $\KK^n$ (for any locally compact field $\KK$ and $n\ge 1$) suffice then to witness such countable equidecomposability \cite[Theorem~1.1]{Tomko_Haar}.

\medskip

Another famous related question is the Banach--Ruziewicz problem, asking whether the Lebesgue measure is the unique (up to normalization) isometry-invariant finitely additive measure defined on bounded Lebesgue measurable subsets of $\RR^n$. The answer is affirmative for $n\ge 3$, which was proved by Margulis \cite{Margulis}; for $n=1,2$ the negative answer is
due to Banach \cite{Banach1923}.
The analogous question was also considered for the spheres
$\mathbb{S}^{n}$, with the $n$-dimensional Hausdorff measure
in place of the Lebesgue measure. The answer is affirmative for $n\ge 2$ due to Margulis \cite{Margulis1980}, Sullivan \cite{Sullivan} and Drinfel'd \cite{Drinfeld}. Later, Bekka established a related result about uniqueness of the Haar measure as a finitely additive $G$-invariant measure for some compact $p$-adic groups $G$ like $\SL(n,\ZZ_p)$, cf. \cite[Theorem~A]{Bekka}.

The positive solution to the Banach--Ruziewicz problem for
$\mathbb{S}^{n}$, $n\ge 2$, is closely connected to some spectral properties (the {\em spectral gap}) of the averaging operator corresponding to the action of a suitable subgroup of $\SO(n+1,\RR)$.
Boutonnet, Ioana and Salehi Golsefidy \cite{BouIoaGol} introduced the notion of the {\em local spectral gap} for general measure-preserving actions. They showed that if $\Gamma$ is a countable dense subgroup of a locally compact second countable group $G$ with the Haar measure $m_G$, then the action $\Gamma\acts (G,m_G)$ has the local spectral gap property if and only if $m_G$ is unique among finitely additive $\Gamma$-invariant measures defined on bounded measurable subsets of $G$ \cite[Theorem~D]{BouIoaGol}. More recently, Salehi Golsefidy provided a general condition under which a $p$-adic analytic group has the local spectral gap property \cite[Theorem~3]{Golsefidy_p-adic}.

The authors of \cite{GrMaPikhurko} established the equivalence of their expansion property and the local spectral gap property \cite[Lemma~5.2]{GrMaPikhurko}. They also showed that, under some mild assumptions, if a measure-preserving action $G\acts (X,\mu)$ is both paradoxical and expanding, then $\mu$ is unique among all finitely additive $G$-invariant measures defined on $\mathcal{D}$, the family of $\mu$-measurable subsets of $X$ that have the Baire property and compact closure \cite[Lemma~8.1]{GrMaPikhurko}. Their proof also works if the family $\mathcal{D}$ is enlarged by dropping the Baire property requirement. Informally, the Banach--Ruziewicz problem for the action $G\acts X$ has a positive answer.
Applying this machinery, our Theorem~\ref{th:main_expanding} can be used to solve the ``$p$-adic'' Banach--Ruziewicz problem as follows (see Corollary~\ref{cor:final}~(ii)).

\begin{namedtheorem}\label{thm:B}
    Let $\pp$ be a non-Archimedean place of a number field $K$ and $n\ge 2$. Let $\Gamma:=\SL(n,\OO_K)\ltimes K^{n}$, where $\OO_K$ is the ring of integers of $K$. The Haar measure on $K_{\pp}^n$ is (up to a constant factor) the unique nonzero, finitely additive, $\Gamma$-invariant, $[0,\infty)$-valued function defined on the family of Haar measurable subsets of $K_{\pp}^n$ with compact closure.
\end{namedtheorem}

We also address the one-dimensional setting: for a locally compact non-Archimedean valued field $\KK$ (of arbitrary characteristic)
whose residue field is not a prime field, we obtain analogous results on measurable
equidecomposability and the Banach--Ruziewicz problem for $\KK$ itself, with respect to
a countable group of not necessarily affine (or even additive) isometries (Corollary~\ref{cor:general_isometries}). The proof uses some properties of unramified extensions of $\QQ_p$ and the
isometric characterization of locally compact non-Archimedean valued fields.

\medskip

The paper is organized as follows. In Section~2.1, we introduce basic concepts of non-Archimedean analysis, focusing mainly on finite-dimensional normed spaces over locally compact fields.
This material is standard; for details we refer to e.g. \cite{UC} and \cite{Neukirch}.
Section~2.2 contains definitions and basic facts on (measurable) equidecomposability, paradoxical and expanding actions, and (local) spectral gap. We also prove some auxiliary results (Lemmas~\ref{lem:equi_from_paradoxicality}, \ref{lem:cover_each_other_dense} and Proposition~\ref{prop:paradoxical_action_of_SA_dense}). Then we recall Theorems~\ref{thm:Golsefidy's_theorem} (due to Salehi Golsefidy) and \ref{thm:consequences_of_parad+expanding} (due to Grabowski, M\'ath\'e and Pikhurko), which are deep results on which we rely later.
Section~3 begins with Proposition~\ref{prop:SG_for_group_acting_on_itself}, establishing the spectral gap for the action $\SL(n,\OO_{K})\acts \SL(n,\OO_{K_{\pp},\pp})$, $n\ge 2$. Then we prove Proposition~\ref{prop:transfer_of_expansion}, which states that, under some general conditions, an expanding action of a subgroup $\Gamma\le G$ on $G$ gives rise to an expanding action of $\Gamma$ on a space acted upon by $G$. We then prove our main Theorem~\ref{th:main_expanding}, and derive from it---as Corollary~\ref{cor:final}---the above mentioned Theorems \ref{thm:A} and \ref{thm:B}.
In Corollaries \ref{cor:unramified} and \ref{cor:general_isometries}, we address analogous problems in the one-dimensional setting.
Finally, we pose two open questions.

\section{Preliminaries}

\subsection{Non-Archimedean analysis}

A {\em non-Archimedean valuation} on a field $\KK$ is any function $\lvert\blank\rvert\colon \KK\to[0,\infty)$ satisfying for all $x,y\in \KK$ the following conditions:\\
\begin{enumerate*}[label=(\arabic*)]
    \item $\lvert x \rvert=0 \Leftrightarrow x=0$,
    \item $\lvert xy \rvert=\lvert x \rvert\lvert y \rvert$,
    \item $\lvert x+y \rvert\leq \max\{\lvert x \rvert,\lvert y \rvert\}$.
\end{enumerate*}\\
Throughout, by a {\em non-Archimedean valued field} we mean a field $\KK$ with a non-Archimedean valuation $\lvert\blank\rvert$ that is {\em nontrivial}, i.e., $\lvert x\rvert \not\in \{0,1\}$ for some $x\in\KK$. The map $(x,y)\mapsto  \lvert x-y \rvert$ is an ultrametric on $\KK$; it turns $\KK$ into a topological field. The multiplicative group $\lvert\KK^{\times}\rvert:=\{\lvert x\rvert: x\in \KK \setminus\{0\}\}$ is called the {\em value group} of $\KK$.
The metric completion of $\KK$ is in a natural way a non-Archimedean valued field.
Two non-Archimedean valuations $\lvert\blank\rvert_1$ and $\lvert\blank\rvert_2$ are called {\em equivalent} if they induce the same topology on $\KK$, which holds precisely if $\lvert\blank\rvert_1=\lvert\blank\rvert_2^s$ for some real number $s>0$.
A {\em non-Archimedean place} of $\KK$ is an equivalence class of nontrivial non-Archimedean valuations on $\KK$.

Let $\pp$ be a non-Archimedean place of a field $\KK$, and let $\lvert\blank\rvert_{\pp}$ be any representative of $\pp$. The closed unit ball $\OO_{\KK,\pp}:=\{x\in \KK: \lvert x\rvert_{\pp} \le 1\}$ forms a subring of $\KK$ called the {\em valuation ring}; it is a local ring with the open unit ball $\mm_{\KK,\pp}:=\{x\in \KK: \lvert x\rvert_{\pp} < 1\}$ as its unique maximal ideal. The unit sphere $\OO_{\KK,\pp}\setminus \mm_{\KK,\pp}$ forms the group of units of $\OO_{\KK,\pp}$. The quotient ring $\OO_{\KK,\pp}/\mm_{\KK,\pp}$ is called the {\em residue field} of $\KK$ (with respect to $\pp$).
It is known that $\KK$ is locally compact if and only if it is complete with respect to $\lvert\blank\rvert_{\pp}$, its residue field is finite, and the value group $\lvert\KK^{\times}\rvert_{\pp}$ is infinite cyclic {\cite [Corollary~12.2]{UC}}.

By a {\em number field} $K$ we mean any field extension $\QQ\subseteq K$ of finite degree. The {\em ring of integers} $\OO_{K}$ of a number field $K$ consists of all $x\in K$ that are integral over $\ZZ$, i.e., such that $f(x)=0$ for a monic polynomial $f\in\ZZ[X]$. There is a well-known bijective correspondence between non-Archimedean places of a number field $K$ and maximal ideals of $\OO_K$, namely: $\pp \mapsto \mm_{K,\pp} \cap \OO_{K}$. Moreover, $\OO_{K}$ is equal to the intersection of $\OO_{K,\pp}$ over all non-Archimedean places $\pp$ of $K$.
If $\pp$ is a non-Archimedean place of $K$, we write $K_{\pp}$ for the completion of $K$ with respect to a fixed representative $\lvert\blank\rvert_{\pp}$, and use the symbols $\lvert\blank\rvert_{\pp}$ and $\pp$ also for the extended valuation on $K_{\pp}$ and the corresponding place.
The valued field $(K_{\pp},\lvert\blank\rvert_{\pp})$ is locally compact and every bounded subset of $K_{\pp}$ has compact closure.
The subring $\OO_{K}$ is dense in $\OO_{K_{\pp},\pp}$ with the topology induced by $\lvert\blank\rvert_{\pp}$. Indeed, it is enough to show that $\OO_K$ is dense in $\OO_{K,\pp}$, which follows from the surjectivity of the natural homomorphism $\OO_K \to \OO_{K,\pp}/\mm_{K,\pp}^n$ for all $n \ge 1$, cf. \cite[Corollary~I.11.2]{Neukirch}. 

Consider the field $K:=\QQ$. Then $\OO_{\QQ}=\ZZ$, and any non-Archimedean place of $\QQ$ is determined by a unique prime number $p$ via the {\em $p$-adic valuation $\lvert\blank\rvert_p$}, defined as follows:
\begin{equation*}
    \lvert 0 \rvert_p:= 0, \; \left\lvert p^k \frac{m}{n}\right \rvert_p:=p^{-k} \quad \text{for } k,m,n\in\mathbb{Z}, \quad p\nmid m, p\nmid n.
\end{equation*}
The completion $\QQ_p$ is called the {\em field of $p$-adic numbers}; we also write $\ZZ_p:=\OO_{\QQ_p,p}$.

Let $\FF_q$ be the finite field with $q$ elements. The {\em field of Laurent series} $\FF_q((T))$ over $\FF_q$ is the field of formal sums $x=\sum_{n\in\ZZ} x_n T^n$, $x_n\in\FF_q$, with well-ordered support $\supp{x}:=\{n\in\ZZ: x_n\neq 0\}$. For any number $c\in(0,1)$,
we introduce the following non-Archimedean valuation on $\FF_q((T))$:
\begin{equation*}
  \left\lvert \sum_{n\in \ZZ}x_n T^n
 \right\rvert_c:= c^{\min \supp{x}} \, \text{ if } x\ne 0, \quad \lvert 0\rvert_c:=0.  
\end{equation*}

Let $(\KK,\lvert\blank\rvert_{\pp})$ be a complete non-Archimedean valued field. Then $\lvert\blank\rvert_{\pp}$ extends uniquely to any algebraic extension $\LL$ of $\KK$ \cite[Theorem~14.2]{UC}.
A finite extension $\KK\subseteq \LL$ of degree $n$ is called {\em unramified} if the extension of residue fields $\OO_{\KK,\pp}/\mm_{\KK,\pp} \subseteq \OO_{\LL,\pp}/\mm_{\LL,\pp}$ has degree $n$ and is separable (separability is automatic if $\OO_{\KK,\pp}/\mm_{\KK,\pp}$ is finite); then $\lvert \KK^{\times}\rvert_{\pp}=\lvert \LL^{\times}\rvert_{\pp}$. For every prime $p$ and $n\ge 1$, adjoining a primitive $(p^n-1)$-th root of unity $\zeta$ to $\QQ_p$ yields an unramified extension $\QQ_p(\zeta)$ of degree $n$ with residue field $\FF_{p^n}$ \cite[Proposition~II.7.12]{Neukirch}.

Recall the following characterization of locally compact non-Archimedean valued fields $(\KK,\lvert\blank\rvert)$ \cite[Proposition~II.5.2]{Neukirch}.
If $\KK$ is of characteristic $0$, then it is a
finite extension of $\QQ_p$ for some prime $p$, and $\lvert\blank\rvert$ is, up
to equivalence, the unique extension of $\lvert\blank\rvert_p$ to $\KK$.
Moreover, $\KK$ contains a dense number field $K$ such that $[K:\QQ]=[\KK:\QQ_p]$, see \cite[Ch.~II, \S~2, Corollary on p.~44]{LangANT},
and hence $\KK=K_{\pp}$ for the non-Archimedean place $\pp$ of $K$ induced by
$\lvert\blank\rvert$.
On the other hand, if $\KK$ has positive characteristic, then it is isometrically isomorphic to $(\FF_q((T)),\lvert\blank\rvert_c)$ for a suitable $c$ and $\FF_q$ equal to the residue field of $\KK$.

We will use the fact that two locally compact non-Archimedean valued fields are isometric if and only if they have equal value groups and their residue fields have the same cardinality; a suitable isometry (nonadditive, in general) can be constructed using a well-known series expansion \cite[Theorem~12.1]{UC}, see also \cite[Corollary~3]{KO-CMB}. Note that passing from $\lvert\blank\rvert$ to an equivalent valuation $\lvert\blank\rvert^s$ on $\KK$, $s>0$, does not change the isometry group.

\medskip

Let $X$ be a linear space over a non-Archimedean valued field $(\KK,\lvert\blank\rvert)$. A function $\lVert \blank \rVert \colon X\to [0, \infty)$ is called a {\em non-Archimedean norm} on $X$ if, for all $x,y \in X$ and $\alpha\in \KK$, we have the following conditions:\\
\begin{enumerate*}[label=(\arabic*)]
   \item $\lVert x \rVert =0 \Leftrightarrow x=0$,
   \item $\lVert \alpha x \rVert = \lvert \alpha \rvert \lVert x \rVert $,
   \item $\lVert x+y \rVert  \leq \max\{\lVert x \rVert ,\lVert y \rVert \}$.
\end{enumerate*}\\
The norm $\lVert \blank \rVert $ induces an ultrametric on $X$ via the map $(x,y)\mapsto \lVert x-y \rVert $.

For any $n\in\NN$, we consider the following {\em maximum norm} on the linear space $\KK^{n}$: $\lVert x \rVert_{\infty} :=\max \{\lvert x_1 \rvert, \dots, \lvert x_n \rvert\}$ for $x=(x_1,\dots,x_n)\in \KK^n$. If $(\KK,\lvert\blank\rvert)$ is locally compact, the resulting topological space $\KK^n$ is locally compact, and every bounded subset of $\KK^n$ has compact closure.
Recall that, in every ultrametric space, all balls (closed or open) and all spheres are both closed and open in the induced topology.

\medskip
If $R$ is a commutative ring with identity and $n\in\NN$, we denote as usual by $\GL(n,R)$ ($\SL(n,R)$) the general (special) linear group of degree $n$ over $R$. 
Let $\KK$ be a field and $n\in \mathbb{N}$.
Recall that an {\em affine transformation} of $\KK^n$ is a  mapping $g\colon \KK^n \to \KK^n$ of the form
$\KK^n\ni x\mapsto L(g)x+\tau(g)$,
where $L(g)\in \M(n,\KK)$ is a linear mapping and $\tau(g) \in \KK^n$ is a fixed vector.
Note that $g$ is invertible if and only if $L(g)\in\GL(n,\KK)$. The group of all invertible affine transformations of $\KK^n$ is isomorphic, via the map $g\mapsto (L(g),\tau(g))$, to the semidirect product $\GL(n,\KK)\ltimes \KK^n$ with respect to the natural action of $\GL(n,\KK)$ on the additive group $\KK^n$.
If $H_1$ is a subgroup of $\GL(n,\KK)$ and $H_2$ is an $H_1$-invariant subgroup of $\KK^n$, then $H_1\ltimes H_2$ is a subgroup of $\GL(n,\KK)\ltimes \KK^n$.

Let $(\KK,\lvert\blank\rvert_{\pp})$ be a locally compact non-Archimedean valued field. For any $n\in\NN$, we give $\M(n,\KK)\cong\KK^{n^2}$ the product topology and $\GL(n,\KK)$ the subspace topology. We also give $G:=\GL(n,\KK)\ltimes\KK^n$ the product topology. Then $G$ is a locally compact, $\sigma$-compact Polish topological group, and its subgroup $\SL(n,\OO_{\KK,\pp})$ is compact.
If $R\subseteq \OO_{\KK,\pp}$ is a dense subring, then $\SL(n,R)$ is a dense subgroup of $\SL(n,\OO_{\KK,\pp})$ (e.g., it follows from the fact that the group $\SL(n,\OO_{\KK,\pp})$ is generated by elementary matrices \cite[Theorem~4.3.9]{HahnOMeara}).

We recall an important characterization of the linear and affine isometries of the normed space $(\KK^n,\lVert\blank\rVert_{\infty})$ for $n\in\NN$. Namely, a linear endomorphism is an isometry of $(\KK^n,\lVert\blank\rVert_{\infty})$ if and only if it is represented (in the canonical basis) by a matrix from $\GL(n,\OO_{\KK,\pp})$, see \cite[Theorem~2.2]{KO-RACSAM}. Consequently, $\GL(n,\OO_{\KK,\pp})\ltimes \KK^n$ can be viewed as the group of affine isometries of $\KK^n$.

\subsection{Equidecomposability and measure theory}

By a {\em group action} $G\acts X$ of a group $G$ on a set $X$ we mean a function $G\times X \ni (g,x)\mapsto g.x\in X$ such that $(gh).x=g.(h.x)$ and $1.x=x$ for all $g,h\in G$, $x\in X$. We put $g.U:=\{g.x: x\in U\}$ and $S.U:=\bigcup_{g\in S}g.U$ for any $U\subseteq X$ and $S\subseteq G$. The action is called {\em transitive} if $G.\{x\}=X$ for some (equivalently, for all) $x\in X$. We say that $A\subseteq X$ {\em $G$-covers} $B\subseteq X$ if there is a finite set $S\subseteq G$ such that $B\subseteq S.A$.

An {\em isomorphism} of group actions $G_1\acts X_1$ and $G_2 \acts X_2$ is a pair $(f,\Phi)$ consisting of a bijection $f\colon X_1 \to X_2$ and a group isomorphism $\Phi\colon G_1\to G_2$ such that
$\Phi(g).f(x)=f(g.x)$ for all $g\in G_1$, $x\in X_1$.

\begin{definition}[{\cite[Definition 3.4 and Corollary 3.7]{W}}]\label{def:equidecomp}
   Let $G\acts X$ be a group action. We say that two sets $A, B\subseteq X$ are {\em $G$-equidecomposable} if, for some $n\in\mathbb{N}$, there exist partitions
\begin{equation}\label{eq:equidecomp}
        A=\bigsqcup_{i=1}^n A_i, \quad B=\bigsqcup_{i=1}^n B_i,
\end{equation}
and elements $g_1,\dots,g_n\in G$ such that $g_i.A_i=B_i$ for all $1\leq i\leq n.$ Then we write $A\sim_{G} B$, and thus obtain an equivalence relation $\sim_{G}$ on the power set of $X$. 

A set $E\subseteq X$ is {\em $G$-paradoxical} if $E=A\sqcup B$ for some subsets $A,B\subseteq E$ with $A\sim_G E \sim_G B$.
\end{definition}

It is known that if a set $E$ is $G$-paradoxical, $A\subseteq E \subseteq X$, and $A$ $G$-covers $E$, then $A\sim_{G} E$ \cite[Corollary~10.22]{W}. We note the following generalization, where the inclusion $A\subseteq E$ is replaced by the weaker hypothesis that $E$ $G$-covers $A$.

\begin{lemma}\label{lem:equi_from_paradoxicality}
 Let $G\acts X$ be a group action, and let $A$, $E$ be two subsets of $X$ that $G$-cover each other. If $E$ is $G$-paradoxical, then $A\sim_{G}E$.
\end{lemma}

\begin{proof}
We adapt the proof of \cite[Corollary~10.22]{W}, with one necessary modification. Since $E$ $G$-covers $A$ and $E$ is $G$-paradoxical, we have $[A]\le m[E]=[E]$ for some $m\in \NN$ in the type semigroup $\mathcal{S}$, and hence $[A]\le [E]$ (in the cited argument, this last inequality followed instead from $A\subseteq E$). All other steps can be kept intact.
\end{proof}

The following definition is motivated by \cite[Definition 1.7]{GrMaPikhurko}, and refers to the phenomenon known also as the ``strong'' Banach--Tarski paradox.

\begin{definition}\label{def:paradoxical_action}
    Let $G\acts X$ be a group action on a locally compact Hausdorff space $X$ by homeomorphisms. The action is called {\em paradoxical} if any two subsets $A,B\subseteq X$, each with compact closure and nonempty interior, are $G$-equidecomposable.
\end{definition}

In this terminology, the famous result of Banach and Tarski \cite[Theorem~24]{BT} says that the action $\SO(n,\RR)\ltimes \RR^n\acts \RR^n$ of the group of orientation-preserving isometries of the Euclidean space $\RR^n$ is paradoxical for $n\ge 3$. For other examples of paradoxical actions, see \cite[Example~1.9]{GrMaPikhurko}.
Note that if $(f,\Phi)$ is an isomorphism of group actions $G_i\acts X_i$ on topological spaces $X_i$ ($i=1,2$) and $f\colon X_1\to X_2$ is a homeomorphism, then $G_1\acts X_1$ is paradoxical if and only if so is $G_2\acts X_2$.

Given a group action $G\acts X$ on a topological space $X$, we denote by $\mathcal{C}=\mathcal{C}(G\acts X)$ the family of all subsets $C\subseteq X$ that have compact closure and $G$-cover a nonempty open set. If $G$ is a topological group, the action $G\acts X$ is called {\em continuous} if it is continuous as a function $G\times X\to X$.

\medskip

We will often use the following condition, under which any two members of $\mathcal{C}(\Gamma \acts X)$ $\Gamma$-cover each other for a dense subgroup $\Gamma\le G$.

\begin{lemma}\label{lem:cover_each_other_dense}
    Let $G\acts X$ be a transitive, continuous action of a topological group $G$ on a topological space $X$. Let $\Gamma$ be a dense subgroup of $G$, and consider the restricted action $\Gamma\acts X$.
    Then, any two sets $A,B\in \mathcal{C}(\Gamma \acts X)$ $\Gamma$-cover each other.
\end{lemma}

\begin{proof}
Take any $A,B\in \mathcal{C}(\Gamma \acts X)$. By definition, $B$ $\Gamma$-covers a nonempty open set $U\subseteq X$. It suffices to show that $U$ $\Gamma$-covers $A$.
Note that, for each $x\in X$, the map $G\ni g\mapsto g^{-1}.x\in X$ is a continuous surjection, under which the preimage of $U$ is nonempty and open, so it meets $\Gamma$. Hence, $X=\Gamma.U$ and the claim follows from the compactness of $\bar{A}$.
\end{proof}

The following observation is a consequence of Lemmas \ref{lem:equi_from_paradoxicality} and \ref{lem:cover_each_other_dense}.

\begin{corollary}\label{cor:paradoxical_action_from_one_paradoxical_set}
    Let $G\acts X$ be a transitive, continuous action of a topological group $G$ on a locally compact Hausdorff space $X$; and let $\Gamma$ be a dense subgroup of $G$. If there exists a $\Gamma$-paradoxical set $E\in\mathcal{C}(\Gamma \acts X)$, then the action $\Gamma\acts X$ is paradoxical.
\end{corollary}

In our paper \cite{KO-RACSAM}, we proved that the affine action of the group $\SL(n,\OO_{\KK,\pp}) \ltimes \KK^n$, denoted there by $\SA(n,D,\KK)$, on $\KK^n$ is paradoxical for any locally compact non-Archimedean valued field $(\KK,\lvert\blank\rvert_{\pp})$ and  $n\ge 2$ \cite[Theorem~5.7]{KO-RACSAM}. In particular, this applies to $\KK:=K_{\pp}$ for a non-Archimedean place $\pp$ of a number field $K$. Let us strengthen this result, replacing $\SL(n,\OO_{\KK,\pp}) \ltimes \KK^n$ with its dense subgroup $\SL(n,\OO_{K})\ltimes K^n$.

\begin{proposition}\label{prop:paradoxical_action_of_SA_dense}
    Let $\pp$ be a non-Archimedean place of a number field $K$.
    The affine action $\SL(n,\OO_K)\ltimes K^n\acts K_{\pp}^n$ is paradoxical for every $n\ge 2$.
\end{proposition}

\begin{proof}
    Let $n\ge 2$, and denote $X:=K_{\pp}^n$, $G:=\SL(n,\OO_{K_{\pp},\pp}) \ltimes K_{\pp}^n$, $\Gamma:=\SL(n,\OO_K)\ltimes K^n$. Since $G$ contains all translations of $K_{\pp}^n$, the action $G\acts X$ is transitive. The ring of integers $\OO_K$ is dense in $\OO_{K_{\pp},\pp}$ and $K$ is dense in $K_{\pp}$, so $\Gamma$ is a dense subgroup of $G$.
    
    Let $E\subseteq X$ be any open ball in $K_{\pp}^n$ not containing the zero vector. By \cite[Theorem~5.1]{KO-RACSAM}, $E$ is $\SL(n,\ZZ)$-paradoxical, and hence $\Gamma$-paradoxical. Since $E\in\mathcal{C}(\Gamma \acts X)$, Corollary~\ref{cor:paradoxical_action_from_one_paradoxical_set} implies that the action $\Gamma\acts X$ is paradoxical.
\end{proof}

\bigskip

A group action $G\acts (X,\mathcal{A})$ on a measurable space is called {\em measurable} if $g.A\in \mathcal{A}$ for all $g\in G$ and $A\in\mathcal{A}$. An action $G\acts (X,\mathcal{A},\mu)$ on a measure space is called {\em measure-preserving} if it is measurable and $\mu(g.A)=\mu(A)$ for all $g\in G$ and $A\in\mathcal{A}$. We will write $\mathcal{A}_{\mu}$ for the completion of $\mathcal{A}$ with respect to $\mu$; the unique extension of $\mu$ onto $\mathcal{A}_{\mu}$ will also be denoted by $\mu$. Note that any measure-preserving action $G\acts(X,\mathcal{A},\mu)$ can be regarded as an action $G\acts(X,\mathcal{A}_{\mu},\mu)$; so we will often write briefly $G\acts (X,\mu)$. Members of $\mathcal{A}_{\mu}$ will be called {\em $\mu$-measurable} sets.

\begin{definition}\label{def:measurable_equidecom}
    Let $G\acts (X,\mathcal{A})$ be a measurable group action. We say that two sets $A, B\in \mathcal{A}$ are {\em $\mathcal{A}$-measurably $G$-equidecomposable} if there are partitions \eqref{eq:equidecomp} and elements $g_1,\dots,g_n\in G$ such that $A_i\in\mathcal{A}$ and $g_i.A_i=B_i$ for all $1\leq i\leq n$.
\end{definition}

If $X$ is a topological space, we write $\mathcal{B}=\mathcal{B}(X)$ for its Borel $\sigma$-algebra, and $\mathcal{BP}$ for the $\sigma$-algebra of all subsets $A\subseteq X$ with the {\em Baire property} (i.e., of the form $A=U\,\triangle\, M$, where $U\subseteq X$ is open and $M$ is a countable union of nowhere dense subsets of $X$). 

\medskip

We recall the important notion of a domain of expansion. It was introduced by Grabowski, M\'ath\'e and Pikhurko \cite[Equation (4)]{GrMaPikhurko}. They used the condition (i) below in their definition. More recently, Li, Vigolo and Zhang used the condition (ii) in \cite[Definition 2.6]{LiVigZhang}, and pointed out that the two definitions are equivalent \cite[Remark 2.7]{LiVigZhang}.

\begin{definition}\label{def:domain_of_exp}
    Let $G\acts (X,\mathcal{A},\mu)$ be a measure-preserving group action.
    A set $C\in \mathcal{A}_{\mu}$ is a {\em domain of expansion} if $0<\mu(C)<\infty$ and one of the following equivalent conditions holds:
    \begin{enumerate}[label=\textup{(\roman*)}]
        \item For every real number $\eta>0$, there exists a finite set $R\subseteq G$ such that, for all $Y\subseteq C$, $Y\in \mathcal{A}_\mu$, we have
        \[
        \mu((R.Y)\cap C) \ge \min\left\{(1-\eta)\mu(C), \frac{\mu(Y)}{\eta}\right\}.
        \]
        \item There exist a real number $c>0$ and a finite set $S\subseteq G$ such that, for all $Y\subseteq C$, $Y\in \mathcal{A}_\mu$, with $\mu(Y)\le \frac{\mu(C)}{2}$, we have
        \begin{equation}\label{eq:domain_of_exp_2nd_cond}
            \mu((S.Y)\cap C)\ge (1+c)\mu(Y).
        \end{equation}
    \end{enumerate}
\end{definition}

\begin{remark}\label{rem:domain_only_borel_sets}
    It is easy to verify that it suffices to check the condition \eqref{eq:domain_of_exp_2nd_cond} only for sets $Y\in \mathcal{A}$.
\end{remark}

The following definition is taken from \cite[Section~1.4]{GrMaPikhurko}.

\begin{definition}\label{def:expanding_action}
    Let $\mu$ be a $\sigma$-finite Borel measure on a Polish topological space $X$. Assume that $\mu$ is nonzero on all nonempty open sets. A measure-preserving group action $G\acts (X,\mu)$ is called {\em expanding} if every $C\in \mathcal{B}_{\mu}\cap \mathcal{C}$ is a domain of expansion.
\end{definition}

\begin{remark}\label{rem:expanding_one_set}
    Assume additionally that the action in Definition \ref{def:expanding_action} is such that any two sets $A,B\in\mathcal{C}$ $G$-cover each other.
    Then, according to \cite[Lemma~2.5]{GrMaPikhurko}, it is enough to check that at least one member of $\mathcal{B}_{\mu}\cap \mathcal{C}$ is a domain of expansion to ensure that the action in question is expanding.
\end{remark}

One can routinely check that the expansion property is stable under certain isomorphisms of group actions in the following sense.

\begin{remark}\label{rem:expanding_stable_under_isomorphisms}
    Let $G_1\acts (X_1,\mu_1)$ be an expanding action, and let $X_2$ be a topological space. Suppose that $(f,\Phi)$ is an isomorphism of group actions $G_1\acts X_1$ and $G_2\acts X_2$ with $f\colon X_1\to X_2$ being a homeomorphism. Denote by $\mu_2:=f_{*}\mu_1$ the {\em push-forward} of $\mu_1$ under $f$, i.e., $\mu_2(A):=\mu_1(f^{-1}(A))$ for all Borel sets $A\subseteq X_2$. Then the action $G_2\acts (X_2,\mu_2)$ is expanding.
\end{remark}

We recall the concepts of the {\em spectral gap} and the {\em local spectral gap}. The former notion is classical, whereas the latter was introduced in \cite{BouIoaGol} and can be regarded as a localized version of the spectral gap. For convenience, we collect them both in one definition.

\begin{definition}[{\cite[Definitions 1.1 and 1.2]{BouIoaGol}}]\label{def:local_spectral_gap}
    Let $G\acts (X,\mathcal{A},\mu)$ be a measure-preserving group action and let $C\in\mathcal{A}_\mu$ with $0<\mu(C)<\infty$. We say that the action has {\em local spectral gap} with respect to $C$ if there exist a finite set $S\subseteq G$ and a real number $\kappa>0$ such that
    \begin{equation}\label{eq:local_spectral_gap}
    \lVert f \rVert_{2,C}\le \kappa \sum_{g\in S}\lVert g. f - f\rVert_{2,C} \quad \text{for all } f\in L^2(X,\mu)  \text{ with } \int_{C}f\, d\mu = 0,
    \end{equation}
    where $(g.f)(x):=f(g^{-1}.x)$ for $g\in G$, $x\in X$, and $\lVert\varphi\rVert_{2,C}$ denotes the $L^2$\nobreakdash-norm of the restriction of $\varphi\in L^2(X,\mu)$ to $C$.

    If $\mu(X)<\infty$, the action is said to have {\em spectral gap} if it has local spectral gap with respect to $X$.
\end{definition}

\begin{remark}\label{rem:operator_norm_SG}
    Let $0<\mu(X)<\infty$ and $\emptyset\ne S\subseteq G$ be a finite symmetric set (i.e., $S^{-1}=S$). Consider the self-adjoint {\em averaging operator} $T_S$ on $L^2(X,\mu)$ given by
    $T_S (f):=\frac{1}{\lvert S\rvert}\sum_{g\in S}g.f$. It has operator norm $1$ and the constant function $\mathbf{1}_X$ is its eigenfunction with eigenvalue $1$.
    Let also $L^{2}_{0}(X,\mu):=\{f\in L^2(X,\mu): \int_{X}f \,d\mu =0\}$. Note that the subspace $L^2_0(X,\mu)$ is invariant under $T_S$. 
   As pointed out in \cite[p.~717]{BouIoaGol}, the existence of $\kappa$ in \eqref{eq:local_spectral_gap} for $C=X$ is equivalent to the presence of a ``gap'' right below $1$ in the spectrum of $T_S$. Thus, the action $G\acts X$ has spectral gap if and only if there exists a finite symmetric set $\emptyset \ne S\subseteq G$ with $\lVert T_S\vert_{L^{2}_{0}(X,\mu)}\rVert<1$, cf. \cite[Definition~2.1]{Benoist} or \cite[p.~7191]{Golsefidy_p-adic}.
\end{remark}

According to \cite[\S~7.1]{Folland_RA}, a {\em Radon measure} on a locally compact Hausdorff space $X$ is a Borel measure $\mu$ such that:
\begin{enumerate}[label=\textup{(\roman*)}, nosep]
    \item $\mu(K) < \infty$ for every compact $K \subseteq X$;
    \item $\mu(U) = \sup\{\mu(K) : K \subseteq U,\ K \text{ compact}\}$ for every open $U \subseteq X$;
    \item $\mu(B) = \inf\{\mu(U) : B \subseteq U,\ U \text{ open}\}$ for every Borel $B \subseteq X$.
\end{enumerate}
If $X$ is additionally second countable, then every Borel measure on $X$ satisfying (i) is Radon \cite[Theorem~7.8]{Folland_RA}.

Recall that a (left) Haar measure on a locally compact Hausdorff topological group $G$ is the unique (up to multiplication by a positive constant) nonzero Radon measure $\mu$ on $G$ that is {\em left-invariant}, i.e., $\mu(gA)=\mu(A)$ for all $A\in\mathcal{B}(G)$.

We will use the following ``metric invariance'' property of Haar measure: If $X_1$, $X_2$ are locally compact abelian groups with topologies induced by invariant metrics $d_1$, $d_2$, respectively, and $f\colon X_1 \to X_2$ is a surjective isometry, then the push-forward $f_{*}m_{X_1}$ of a Haar measure on $X_1$ is a Haar measure on $X_2$. It is a consequence of Christensen's theorem \cite[Theorem~1]{Christensen}.
In particular, if $m_n$ is a Haar measure on the additive group $\KK^n$ for a locally compact non-Archimedean valued field $(\KK,\lvert\blank\rvert_{\pp})$, then the (isometric) action $\GL(n,\OO_{\KK,\pp})\ltimes \KK^n\acts (\KK^n,m_n)$ is measure-preserving.

\bigskip

We now cite the main theorem of Salehi Golsefidy's paper \cite{Golsefidy_p-adic}.

\begin{theorem}[{\cite[Theorem 1]{Golsefidy_p-adic}}]\label{thm:Golsefidy's_theorem}
    Let $K$ be a number field, $n\in\NN$, and let $S$ be a nonempty finite symmetric subset of $\GL(n,K)$. Let $\Gamma:=\langle S\rangle$ and $\Gamma_{\pp}$ be the closure of $\Gamma$ in $\GL(n,K_{\pp})$ for any non-Archimedean place $\pp$ of $K$. Denote by $\Sigma$ the set of all $\pp$ such that $\Gamma_{\pp}$ is compact. For any $\pp\in \Sigma$, let $m_{\Gamma_{\pp}}$ be the Haar probability measure on the Borel $\sigma$-algebra of $\Gamma_{\pp}$.
    
    Suppose that the Zariski-closure of $\Gamma$ in $\GL(n,K)$ is a semisimple algebraic group. Then
    $\sup_{\pp\in \Sigma}\lVert T_S\vert_{L^{2}_{0}(\Gamma_{\pp},m_{\Gamma_{\pp}})}\rVert<1$. In particular, the left translation action $\Gamma\acts (\Gamma_{\pp},m_{\Gamma_{\pp}})$ has spectral gap for every $\pp\in \Sigma$.
\end{theorem}

We recall the following result establishing the connection between the notions of a domain of expansion and the local spectral gap. It was proved by Grabowski, M\'ath\'e and Pikhurko, who built upon \cite[Theorem 7.6]{BouIoaGol}.

\begin{lemma}[{\cite[Lemma 5.2]{GrMaPikhurko}}]\label{lem:LSG=domain}
    Let $\mu$ be a $\sigma$-finite Borel measure on a Polish topological space $X$. Assume that $\mu$ is nonzero on all nonempty open sets. Let $G\acts (X,\mu)$ be a measure-preserving group action. 
    
    Then, for any $C\in\mathcal{B}_{\mu}$ with $0<\mu(C)<\infty$, the following are equivalent:
    \begin{enumerate}[label=\textup{(\roman*)}]
        \item The set $C$ is a domain of expansion for the action $G\acts (X,\mu)$.
        \item The action $G\acts (X,\mu)$ has local spectral gap with respect to $C$.
    \end{enumerate}
\end{lemma}

Finally, let us collect two results of \cite{GrMaPikhurko} that highlight the importance of expanding actions in the context of measurable equidecomposability and the Banach--Ruziewicz problem.

\begin{theorem}[{\cite[Corollary 1.10, Lemma 8.1]{GrMaPikhurko}}]\label{thm:consequences_of_parad+expanding}
    Assume that $G\acts (X,\mu)$ is a paradoxical and expanding action. Let $\mathcal{A}\in\{\mathcal{B}_{\mu}, \mathcal{B}_{\mu}\cap \mathcal{BP}\}$.
    \begin{enumerate}[label=\textup{(\roman*)}]
        \item If $A\in \mathcal{A}\cap \mathcal{C}$, then a subset $B\subseteq X$ is $\mathcal{A}$-measurably $G$-equidecomposable to $A$ if and only if $B\in \mathcal{A}\cap \mathcal{C}$ and $\mu(A)=\mu(B)$.
        \item Assume additionally that $\mu$ is atomless and $\mu(U)<\infty$ for some nonempty open set $U\subseteq X$. Let $\mathcal{D}$ be the family of all elements of $\mathcal{A}$ with compact closure. Then every finitely additive $G$-invariant function $\kappa\colon \mathcal{D}\to [0,\infty)$ is a constant multiple of $\mu$.
    \end{enumerate}  
\end{theorem}

We point out that assertion (ii) as stated in \cite[Lemma~8.1]{GrMaPikhurko} treats the $\mathcal{A}=\mathcal{B}_{\mu}\cap \mathcal{BP}$ case. However, the proof can be easily modified to cover the case $\mathcal{A}=\mathcal{B}_{\mu}$ as well.

\section{Results}

We begin with a simple consequence of Theorem~\ref{thm:Golsefidy's_theorem}.

\begin{proposition}\label{prop:SG_for_group_acting_on_itself}
    Let $\pp$ be a non-Archimedean place of a number field $K$, and let $n\ge 2$. The left translation action of $\SL(n,\OO_K)$ on $\SL(n,\OO_{K_{\pp},\pp})$, equipped with the Haar probability measure, has spectral gap.
\end{proposition}

\begin{proof}
Let $\Gamma:=\SL(n,\OO_K)$ and consider the compact group $\Gamma_{\pp}:=\nolinebreak\SL(n,\OO_{K_{\pp},\pp})$ equipped with the Haar probability measure $m_{\Gamma_{\pp}}$ on its Borel $\sigma$-algebra.

It is known that the group $\Gamma$ is finitely generated, cf. \cite{OMeara1965}. Let $S\subseteq \Gamma$ be a finite symmetric set that generates $\Gamma$. Note that $\Gamma_{\pp}$ is the closure of $\Gamma=\left<S\right>$ in $\GL(n,K_{\pp})$.
Moreover, it is an elementary fact that the Zariski-closure of $\Gamma$ in $\GL(n,K)$ equals $\SL(n,K)$, which is known to be a semisimple algebraic group \cite[\S~19.5]{Humph}.

Now, we can directly apply Theorem \ref{thm:Golsefidy's_theorem}, yielding the spectral gap property of the left translation action $\Gamma \acts (\Gamma_{\pp},m_{\Gamma_{\pp}})$.
\end{proof}

Now we prove a result on the ``transfer'' of the expansion property from the left translation action to a general measure-preserving action. The proof is motivated by (and, in fact, generalizes) the argument used in \cite[Section~6.4]{GrMaPikhurko} to prove that the action $\SL(2,\RR)\acts \RR^2\setminus\{0\}$ is expanding.

\begin{proposition}\label{prop:transfer_of_expansion}
    Let $G\acts (X,\mu)$ be a continuous, measure-preserving action of a locally compact, $\sigma$-compact Polish topological group $G$ on a Hausdorff space $X$ with a nonzero Borel measure $\mu$ that is finite on all compact subsets of $X$.
    Let $x_0\in X$ and suppose that there exists a continuous function $s\colon X \to G$ such that $s(x).x_0=x$ for all $x\in X$. Equip $G$ with a left Haar measure $m_G$.

    If, for a dense subgroup $\Gamma\le G$, the left translation action $\Gamma \acts (G,m_G)$ is expanding, then so is the action $\Gamma \acts (X,\mu)$. 
\end{proposition}

\begin{proof}
    Note that the existence of $s$ implies the transitivity of the action $G\acts (X,\mu)$.
    Let $H:=\{g\in G: g.x_0=x_0\}$; it is a closed subgroup of $G$. Equip $H$ with a left Haar measure $m_H$. Define a function $\varphi\colon X\times H \to G$ by
    \[
    \varphi(x,h):=s(x)h \quad \text{for } (x,h)\in X\times H.
    \]
    It can be shown by a direct calculation that $\varphi$ is a homeomorphism with the inverse
    $\varphi^{-1}(g)=(g.x_0, s(g.x_0)^{-1}g)$ for $g\in G$. Since $X$ is homeomorphic to a closed subspace of $G$, it is locally compact, $\sigma$-compact and Polish. The measure $\mu$ is thus $\sigma$-finite. We claim that $\mu(U)>0$ for all nonempty open sets $U\subseteq X$. Suppose that $\emptyset \ne U\subseteq X$ is an open set with $\mu(U)=0$. It follows from the properties of the action $G\acts (X,\mu)$ and the second countability of $X$ that $X=I.U$ for a countable set $I\subseteq G$. Hence, $
    \mu(X)=0$---a contradiction. So far, we have shown that $(X,\mu)$ and $\Gamma\acts (X,\mu)$ satisfy the prerequisites of Definition \ref{def:expanding_action}.

    Notice that $\mu$ is a Radon measure by \cite[Theorem~7.8]{Folland_RA}.
    Since $X$ and $H$ are second countable, the Borel $\sigma$-algebra $\mathcal{B}(X\times H)$ is equal to the product $\sigma$-algebra $\mathcal{B}(X) \otimes \mathcal{B}(H)$, and the product $\mu\times m_H$ of Radon measures $\mu$ and $m_H$ is a Radon measure on $X\times H$ \cite[Theorem~7.20]{Folland_RA}. 
    Let $\mu'$ be the push-forward of $\mu\times m_H$ under the homeomorphism $\varphi$, that is, $\mu'(A):=(\mu\times m_H)(\varphi^{-1}(A))$ for all Borel sets $A\subseteq G$. It is clearly a Radon measure on $G$. We will prove that $m_G=\alpha\,\mu'$ for a constant $\alpha>0$. By the uniqueness of Haar measure, it suffices to show that $\mu'$ is left invariant.

    For any $g\in G$ and $A\in \mathcal{B}(G)$, the equality
    $\mathbf{1}_{\varphi^{-1}(g^{-1}A)}(x,h)=\mathbf{1}_{A}(g\varphi(x,h))=\mathbf{1}_{A}(gs(x)h)$ holds for all $(x,h)\in X\times H$. Hence, by Tonelli's theorem, we have
    \begin{equation}\label{eq:Tonelli_1}
        \begin{split}
            \mu'(g^{-1}A)
            =\int_{X\times H}\mathbf{1}_{\varphi^{-1}(g^{-1}A)}\,d(\mu\times m_H)
            =\int_{X}\int_{H}\mathbf{1}_{A}(gs(x)h)\,dm_H(h)\,d\mu(x).
        \end{split}
    \end{equation}
    It remains to show that the value of \eqref{eq:Tonelli_1} does not depend on $g\in G$.
    We claim that, for any fixed $g\in G$ and $x\in X$, the inner integral in \eqref{eq:Tonelli_1} equals
    $\int_{H}\mathbf{1}_{A}(s(g.x)h)\,dm_H(h)$.
    Indeed, $gs(x).x_0=g.x=s(g.x).x_0$; so the constant $h':=s(g.x)^{-1}gs(x)$  stabilizes $x_0$ and hence belongs to $H$. Therefore,
    \begin{equation}\label{eq:H-invariance}
        \int_{H}\mathbf{1}_{A}(gs(x)h)\,dm_H(h)=\int_{H}\mathbf{1}_{A}(s(g.x)h'h)\,dm_H(h)=\int_{H}\mathbf{1}_{A}(s(g.x)h)\,dm_H(h)
    \end{equation}
    as desired, where the last equality follows from the left $H$-invariance of $m_H$.
     Substituting \eqref{eq:H-invariance} into \eqref{eq:Tonelli_1} and changing the order of integration, we get
    \[\mu'(g^{-1}A)=\int_{H}\int_{X}\mathbf{1}_{A}(s(g.x)h)\,d\mu(x)\,dm_H(h).\]
    Finally, applying the $G$-invariance of $\mu$, we obtain
    \[\mu'(g^{-1}A)=\int_{H}\int_{X}\mathbf{1}_{A}(s(x)h)\,d\mu(x)\,dm_H(h),
    \] which is independent of $g\in G$.

    Without loss of generality, scaling $m_G$ if needed, assume that $\alpha=1$.
    Fix a nonempty open set $C\subseteq X$ with compact closure (it exists by the local compactness of $X$).
    By Lemma \ref{lem:cover_each_other_dense} and Remark \ref{rem:expanding_one_set}, we only need to show that this fixed $C$ is a domain of expansion for the action $\Gamma \acts (X,\mu)$. Clearly, $0<\mu(C)<\infty$.
    Fix also a nonempty open set $V\subseteq H$ with compact closure; note that $0<m_H(V)<\infty$. Define $C':=\varphi(C \times V)\subseteq G$. Since $\Gamma \acts (G,m_G)$ is expanding, the nonempty set $C'$, being open with compact closure, is a domain of expansion for this action. Let $c>0$ and $S\subseteq \Gamma$ be as in Definition \ref{def:domain_of_exp}~(ii) for $C'$ and $\Gamma \acts (G,m_G)$.
    We will prove that the same $c$ and $S$ witness that $C$ is a domain of expansion for the action $\Gamma \acts (X,\mu)$.

    Let $Y\subseteq C$ be a Borel set with $\mu(Y)\le \frac{\mu(C)}{2}$, and put $Y':=\varphi(Y\times V)$.
    Note that $m_G(C')=\mu'(C')=(\mu \times m_H)(C\times V)=\mu(C)m_H(V)$, and similarly $m_G(Y')=\mu(Y)m_H(V)$. Hence, $m_G(Y')\le \frac{m_G(C')}{2}$; and so
    \begin{equation}\label{eq:expansion_for_m_G}
        m_G((SY')\cap C')\ge (1+c)m_G(Y')=(1+c)\mu(Y)m_H(V).
    \end{equation}
    It follows from the formulas expressing $\varphi$ and $\varphi^{-1}$ that $\varphi^{-1}(g\varphi(x,h))=\varphi^{-1}(gs(x)h)\in \{g.x\}\times H$ for all $(x,h)\in X\times H$ and $g\in G$. Therefore, $\varphi^{-1}(SY')\subseteq (S.Y)\times H$, which gives the inclusion
    $\varphi^{-1}((SY')\cap C') \allowbreak\subseteq ((S.Y) \cap\nobreak C)\times V$. We thus obtain
    \[\mu ((S.Y)\cap C) m_H(V) \ge (\mu\times m_H)(\varphi^{-1}((SY')\cap C'))=m_G((SY')\cap C').\] Combining this with $\eqref{eq:expansion_for_m_G}$ and factoring out $m_H(V)$, we get the desired inequality \eqref{eq:domain_of_exp_2nd_cond}.
\end{proof}

We are ready to prove our main theorem.
 
\begin{theorem}\label{th:main_expanding}
    Let $\pp$ be a non-Archimedean place of a number field $K$, and $n\ge 2$. Denote by $X$ the unit sphere in $K_{\pp}^n$, and by $m_n$ the Haar measure on $K_{\pp}^n$. Then both actions $\SL(n,\OO_K)\acts (X,m_n)$, $\SL(n,\OO_K)\ltimes K^{n} \acts (K_{\pp}^{n},m_n)$ are expanding.
\end{theorem}

\begin{proof}
    Note that the actions under consideration preserve $m_n$.
    Let $G:=\SL(n,\OO_{K_{\pp},\pp})$, $\Gamma:=\SL(n,\OO_K)$, and let $m_G$ be the Haar probability measure on $G$. By Proposition \ref{prop:SG_for_group_acting_on_itself}, the left translation action $\Gamma \acts (G,m_G)$ has local spectral gap with respect to $G$. It follows from Lemma \ref{lem:LSG=domain} that $G$ is a domain of expansion for this action. 
    Since the compact set $G$ belongs to the family $\mathcal{C}(\Gamma \acts G)$, the action $\Gamma \acts (G,m_G)$ is expanding by Lemma \ref{lem:cover_each_other_dense} and Remark \ref{rem:expanding_one_set}.

    Put $e_1:=(1,0,\dots,0)$. We will show that there exists a continuous function $s\colon X\to G$ such that $s(x).e_1=x$ for all $x\in X$, which means that $s(x)$ is an $n\times n$ matrix over $\OO_{K_{\pp},\pp}$ with determinant $1$ and first column $x$.

    Define $X_k:=\{(x_1,\dots,x_n)\in X: \min\{1\le i\le n: x_i\notin \mm_{K_{\pp},\pp}\}=k\}$ for $1\le k\le n$. The sets $X_k$ ($1\le k\le n$) form a disjoint partition of $X=\OO_{K_{\pp},\pp}^n\setminus \mm_{K_{\pp},\pp}^n$, and each of them is open as the product $\mm_{K_{\pp},\pp}^{k-1}\times (\OO_{K_{\pp},\pp}\setminus \mm_{K_{\pp},\pp})\times \OO_{K_{\pp},\pp}^{n-k}$. 
    Therefore, it suffices to define $s$ continuously on each $X_k$. Note that, for every $x=(x_1,\dots,x_n)\in X_k$ ($1\le k\le n$), the entry $x_k$ is invertible in $\OO_{K_{\pp},\pp}$.

\medskip
\noindent Case $k=1$. For $x=(x_1,\dots,x_n)\in X_1$, define $s(x)\in \M(n,\OO_{K_{\pp},\pp})$
by
\[
s(x)_{ij}:=
\begin{cases}
    x_i & \text{if } 1\le i\le n, \,j=1;\\
    x_1^{-1} & \text{if } i=j=2;\\
    1 & \text{if } 2<i=j\le n;\\
    0 & \text{otherwise}.
\end{cases}
\]

\medskip
\noindent Case $k>1$. For $x=(x_1,\dots,x_n)\in X_k$ with $k>1$, define
$s(x)\in \M(n,\OO_{K_{\pp},\pp})$ by
\[
s(x)_{ij}:=
\begin{cases}
    x_i & \text{if } 1\le i\le n,\,j=1;\\
    -x_k^{-1} & \text{if } i=1,\, j=k;\\
    1 & \text{if } i=j\notin\{1,k\};\\
    0 & \text{otherwise}.
\end{cases}
\]

It is immediate that $\det{s(x)}=1$ for all $x\in X$. If $k=1$, then $s(x)$ is lower triangular with diagonal entries $x_1,x_1^{-1},1,\dots,1$. If $k>1$,
then the transposition of rows $1$ and $k$ turns $s(x)$ into a lower triangular matrix with
diagonal entries $x_k,1,\dots,1,-x_k^{-1},1,\dots,1$.

It now follows from Proposition \ref{prop:transfer_of_expansion} that the action $\Gamma\acts (X,m_n)$ is expanding, and hence the set $X\in\mathcal{C}(\Gamma\acts X)$ is a domain of expansion for this action. It is clear by Definition \ref{def:domain_of_exp} that $X$ is also a domain of expansion for the action of the larger group $\Gamma':=\Gamma\ltimes K^{n}$ on $(K_{\pp}^{n},m_n)$.

Since the group $G\ltimes K_{\pp}^n$ contains all translations of $K_{\pp}^{n}$, its action on $K_{\pp}^n$ is transitive. Note that $\Gamma'$ is a dense subgroup of $G\ltimes K_{\pp}^n$.
Observe that $X$ is compact and open in $K_{\pp}^{n}$, so it belongs to $\mathcal{C}(\Gamma'\acts K_{\pp}^{n})$. Applying Lemma \ref{lem:cover_each_other_dense} and Remark \ref{rem:expanding_one_set} again, we conclude that the action $\Gamma' \acts (K_{\pp}^{n},m_n)$ is expanding as well.
\end{proof}

Since the action $\SL(n,\OO_K)\ltimes K^{n} \acts (K_{\pp}^{n},m_n)$ is both paradoxical (by Proposition~\ref{prop:paradoxical_action_of_SA_dense}) and expanding (by Theorem~\ref{th:main_expanding}), an application of Theorem~\ref{thm:consequences_of_parad+expanding} yields the following corollary.

\begin{corollary}\label{cor:final}
    Let $\pp$ be a non-Archimedean place of a number field $K$, and let $n\ge 2$. Consider the action $\Gamma\acts (K_{\pp}^{n},m_n)$, where $\Gamma:=\SL(n,\OO_K)\ltimes K^{n}$ and $m_n$ is the Haar measure on $K_{\pp}^{n}$.
    Let $\mathcal{C}:=\mathcal{C}(\Gamma \acts K_{\pp}^n)$ and $\mathcal{A}\in\{\mathcal{B}_{m_n}, \mathcal{B}_{m_n}\cap \mathcal{BP}\}$.
\begin{enumerate}[label=\textup{(\roman*)}]
     \item If $A\in \mathcal{A}\cap \mathcal{C}$, then a subset $B\subseteq K_{\pp}^n$ is $\mathcal{A}$-measurably $\Gamma$-equidecomposable to $A$ if and only if $B\in \mathcal{A}\cap \mathcal{C}$ and $m_n(A)=m_n(B)$.
    \item Let $\mathcal{D}$ be the family of all elements of $\mathcal{A}$ with compact closure. Then every finitely additive $\Gamma$-invariant function $\kappa\colon \mathcal{D}\to [0,\infty)$ is a constant multiple of $m_n$.
\end{enumerate}
\end{corollary}

\begin{example}\label{ex:empty_interior}
    Consider $K:=\QQ$ with the $p$-adic place. Let $A:=\ZZ_p^2$ and $B:=\ZZ_p^2 \setminus \ZZ^2$; they are clearly Borel subsets of $\QQ_p^2$ with compact closure, and $m_2(A)=m_2(B)$. Though $B$ has empty interior, it belongs to $\mathcal{C}(\Gamma \acts \QQ_p^2)$. Indeed, $\ZZ_p^2 \subseteq B \cup g.B$, where $g\in \Gamma$ is the translation by $(1/(p+1),0)\in \QQ^2$. Hence, $A$ and $B$ are $\SL(2,\ZZ)\ltimes \QQ^2$-equidecomposable using Haar measurable pieces with the Baire property. 
\end{example}

Note that every locally compact non-Archimedean valued field $\LL$ of characteristic $0$ can be realized (non-uniquely, in general) as $K_{\pp}$ for a number field $K$ and a non-Archimedean place $\pp$ of $K$. Hence, the conclusions (i), (ii) above hold for $\LL^n$ ($n\ge 2$), with $\Gamma$ a certain countable group of $\LL$-affine isometries of $\LL^n$.

The case $n=1$ is different since the group
$\IA(\LL)$ of $\LL$-affine isometries of $\LL$ is amenable, and hence the action $\IA(\LL)\acts \LL$ is not paradoxical.
We now show that this obstruction
disappears if we consider a field $\LL$ that is a proper unramified extension of a locally compact subfield $\KK$, and $\Gamma$ is a certain group of $\KK$-affine isometries of $\LL$.

\begin{corollary}\label{cor:unramified}
    Let $\KK$ be a locally compact non-Archimedean valued field of characteristic $0$, and let $\KK\subsetneq \LL$ be a proper unramified extension of finite degree. Then there exists a countable group $\Gamma$ of $\KK$-affine isometries of $\LL$ such that the action $\Gamma \acts (\LL,m_{\LL})$, where $m_{\LL}$ is the Haar measure on $\LL$, is paradoxical and expanding. As a consequence, the assertions (i), (ii) of Corollary~\ref{cor:final} hold with $K_{\pp}^n$, $m_n$ replaced by $\LL$ and $m_{\LL}$.
\end{corollary}

\begin{proof}
    Denote by $\lvert\blank\rvert_{\pp}$ the valuation on $\KK$ and its extension to $\LL$.
    Let $n:=[\LL:\KK]\ge 2$. Since the extension $(\KK,\lvert\blank\rvert_{\pp}) \subsetneq (\LL, \lvert\blank\rvert_{\pp})$ is unramified, there exist $b_1,\dots,b_n\in \OO_{\LL,\pp}$ whose images under the residue map $r\colon \OO_{\LL,\pp} \to \OO_{\LL,\pp}/\mm_{\LL,\pp}$ are linearly independent over $\OO_{\KK,\pp}/\mm_{\KK,\pp}$.
    We claim that
    \begin{equation}\label{eq:orthonormal_basis}
        \left\lvert\sum_{i=1}^n \alpha_i b_i \right\rvert_{\pp} = \max_{1\le i\le n} \lvert\alpha_i \rvert_{\pp} \quad \text{for all } \alpha_1,\dots,\alpha_n\in \KK.
    \end{equation}
    Indeed, assume without loss of generality that $\lvert\alpha_1\rvert_{\pp}= \max_{1\le i\le n} \lvert\alpha_i \rvert_{\pp}>0$ and put 
    $x:=\alpha_{1}^{-1}\sum_{i=1}^{n} \alpha_i b_i$.
    Then $x=b_1 + \sum_{i=2}^{n}\beta_i b_i\in \OO_{\LL,\pp}$, where $\beta_i:=\alpha_{1}^{-1}\alpha_i\in \OO_{\KK,\pp}$ for all $2\le i\le n$.
    By $\OO_{\KK,\pp}/\mm_{\KK,\pp}$-linear independence of $r(b_i)$ $(1\le i \le n)$, we obtain $r(x)\ne r(0)$, and so $\lvert x\rvert_{\pp}=1$, which proves \eqref{eq:orthonormal_basis}.

    Hence, the $\KK$-linear isomorphism $f\colon \KK^n\to \LL$, $f(\alpha_1,\dots,\alpha_n):=\sum_{i=1}^n\alpha_i b_i$ is an isometry between $\KK^n$ (with the maximum norm) and $\LL$. 
    Realize $\KK$ as the completion of a number field $K\subseteq \KK$ with respect to $\pp$, and let $\Gamma':=\SL(n,\OO_{K})\ltimes K^n$. For any $g\in \Gamma'$ let $\Phi(g):=f\circ g \circ f^{-1}$, which is a $\KK$-affine isometry of $\LL$. Thus, we obtained an isomorphism $(f,\Phi)$ of group actions $\Gamma'\acts \KK^n$ and $\Gamma \acts \LL$, where $\Gamma:=\Phi(\Gamma')$. The action $\Gamma'\acts (\KK^n,m_n)$ is paradoxical and expanding. By the metric invariance of Haar measure, the push-forward $f_{*}m_n$ is proportional to the Haar measure $m_{\LL}$. Therefore, the action $\Gamma \acts (\LL,m_{\LL})$ is also paradoxical and expanding (see Remark~\ref{rem:expanding_stable_under_isomorphisms}), and the remaining conclusions follow from Theorem~\ref{thm:consequences_of_parad+expanding}.
\end{proof}

Relying on the fact that a locally compact non-Archimedean valued field is, up to isometry, determined by its value group and the cardinality of its residue field, we obtain the following results concerning measurable equidecomposability in such a field (of arbitrary characteristic) with respect to certain of its isometries (which need not be affine or even additive).

\begin{corollary}\label{cor:general_isometries}
    Let $\LL$ be a locally compact non-Archimedean valued field whose residue field is not a prime field. Then there exists a countable group $\Gamma$ of isometries of $\LL$ such that the action $\Gamma \acts (\LL,m_{\LL})$, where $m_{\LL}$ is the Haar measure on $\LL$, is paradoxical and expanding. As a consequence, the assertions (i), (ii) of Corollary~\ref{cor:final} hold with $K_{\pp}^n$, $m_n$ replaced by $\LL$ and $m_{\LL}$.
\end{corollary}

\begin{proof}
    By the assumption, the residue field of $(\LL,\lvert\blank\rvert)$ is isomorphic to $\FF_{p^n}$ for $n\ge 2$ and a prime $p$.
    Applying Corollary~\ref{cor:unramified} to the unramified extension $\QQ_{p}\subsetneq \LL':=\QQ_p (\zeta)$, where $\zeta$ is a primitive $(p^n -1)$-th root of unity, we obtain a countable group $\Gamma'$ of $\QQ_p$-affine isometries of $\LL'$ such that the action $\Gamma' \acts (\LL',m_{\LL'})$ is paradoxical and expanding.

    Let $c>1$ be a generator of $\lvert\LL^{\times}\rvert$ and $s:=\log_{p} c$. Observe that the valued fields $(\LL,\lvert\blank\rvert)$ and $(\LL',\lvert \blank\rvert_p^s)$ have isomorphic residue fields and equal value groups. Hence, there exists a surjective isometry $f\colon \LL' \to \LL$. As in the proof of Corollary~\ref{cor:unramified}, we obtain an appropriate isomorphism $(f,\Phi)$ of group actions $\Gamma'\acts \LL'$ and $\Gamma:=\Phi(\Gamma')\acts \LL$, from which the desired assertions follow.
\end{proof}

\subsection*{Open problems}

Let us suggest the following questions. 

\begin{question}\label{q:positive_characteristic}
    Does an analogue of Theorem~\ref{th:main_expanding} hold for $K$ being a global field of positive characteristic (i.e., a finite extension of the rational function field $\FF_{p}(T)$, $p$ prime)?
\end{question}

\begin{question}\label{q:prime_residue_field}
    Let $\LL$ be a locally compact non-Archimedean valued field with residue field $\FF_p$, $p$ prime. Is the action $\Iso(\LL)\acts (\LL,m_{\LL})$ expanding?
\end{question}

Note that the action in Question~\ref{q:prime_residue_field} is paradoxical \cite[Theorem~11]{KO-CMB}.

\section*{Acknowledgments}

The author would like to warmly thank Oleg Pikhurko for a helpful consultation.

\section*{Statements and Declarations}

{\bf Competing Interests:} The author has no inerests to disclose.

\bibliographystyle{siam}
\bibliography{references}

\end{document}